\documentclass[12pt]{article}

\usepackage{amssymb,amsmath,oldlfont,epsf,stmaryrd}
\usepackage[a4paper,top=3cm,bottom=3cm,left=3cm,right=3cm]{geometry}

\def\fa{\forall}
\def\ex{\exists}
\newcommand{\quot}[1]{\ulcorner #1 \urcorner}
\newcommand{\reflect}[1]{\mbox{`}#1\mbox{'}}

\begin{document}
\title{Teaching G\"odel's incompleteness theorems}
\author{Gilles Dowek\thanks{Inria and 
\'Ecole normale sup\'erieure de Paris-Saclay,
61, avenue du Pr\'esident Wilson,
94235 Cachan Cedex, France,
{\tt gilles.dowek@ens-paris-saclay.fr}.}}
\date{}
\maketitle
\thispagestyle{empty}

\section{Introduction}

The basic notions of logic---predicate logic, Peano arithmetic,
incompleteness theorems, etc.---have for long been an advanced
topic. In the last decades, they became more widely taught, in
philosophy, mathematics, and computer science departments, to graduate
and to undergraduate students. Many textbooks now present these
notions, in particular the incompleteness theorems.

Having taught these notions for several decades, our community can now
stand back and analyze the choices faced when designing such a
course. In this note, we attempt to analyze the choices faced when
teaching the incompleteness theorems. In particular, we attempt to
defend the following points.

\begin{itemize}
\item The incompleteness theorems are a too rich subject to be taught
  in only one course. It is impossible to reach, in a few weeks, the
  second incompleteness theorem and L\"ob's theorem with students who
  have never been exposed to the basics of predicate logic, exactly in
  the same way that it is impossible to reach in a few weeks the
  notion of analytic function with students who have never had been
  exposed to the notions of function and complex number.

  Thus, the incompleteness theorems must be taught several times, at
  different levels, for instance, first, in an elementary course, in
  the third year of university, then in an advanced one, in the fourth
  year. The goals and focus of theses courses are different.

\item When the incompleteness theorems are taught in isolation, they
  are often viewed as a ``promised land'' and the notions of
  computability, representation, reduction, diagonalization, etc. are
  introduced in order to prove these theorems.

  In contrast, we defend that these notions should be motivated and
  taught independently, possibly in different courses. This way, the
  chapter on the incompleteness theorems remains small and focused on
  the specifics of incompleteness.

\item Making the proofs too concrete, for instance by defining
  explicit numberings, so that the students can put their hands on
  this notion, often overloads the proofs with irrelevant
  idiosyncratic details.

  The proofs need to be made abstract enough so that the students can
  focus on the meaningful points.  Tools to make the proof abstract
  enough---abstract syntax, articulation, general fixed-point
  theorems, provability logic notations, etc.---exist and can be used.
\end{itemize}

\section{Which theories ?}

The incompleteness theorems do not apply to all theories. They do not
apply to some theories because they are too weak, such as the theory
containing a binary predicate symbol $=$ and an axiom $\fa x \fa y~(x
= y)$, that is complete. They do not apply to some theories because
they are too strong, such as inconsistent theories, that also are
complete.

\subsection{Strong enough}

Although G\"odel's original proofs applied to the {\it Principia
  Mathematica}, a natural choice for proving the incompleteness
theorems is Peano arithmetic, that already permits to represent
computable functions.

Yet, for the first incompleteness theorem, it is traditional to
consider a much weaker theory: Robinson's arithmetic that is
essentially Peano arithmetic, minus the induction axiom, plus a few
consequences of induction.  This generalization is indeed interesting,
but it should not be considered as mandatory, when teaching the
incompleteness theorem, specially in an elementary course, as it may
divert the students from the meaningful points in the proof of this
theorem.

Moreover, proving the Hilbert-Bernays lemmas and hence the second
incompleteness theorem seems to require a theory stronger than
Robinson's arithmetic, for instance Peano arithmetic.

\subsection{Weak enough}

A natural, but false, idea of the students discovering the incompleteness
of arithmetic is that some axiom has been forgotten and that adding
this axiom will make the theory complete.

So, it is important to remark that adding a finite number of axioms
cannot make arithmetic complete. The right way to formulate this
remark is to prove essential incompleteness, that is that any
extension of arithmetic verifying, some properties, is incomplete.
The first condition can be formulated as the fact that the set of
axioms is decidable, that proof-checking is decidable, or as the fact
that the set of theorems is semi-decidable. The second is the
consistency of the theory.
 
Yet, assuming only consistency requires to prove the G\"odel-Rosser
theorem, introducing some unneeded complexity. An alternative is to
assume $\omega$-consistency, that is that each time a proposition of
the form $\ex x~A$ is provable, there exists a natural number $n$ such
that $\neg (\underline{n}/x)A$ is not provable. But, this property
also introduces some unneeded complexity.

A simpler option is to assume a stronger hypothesis: that the theory has a 
standard model.
An advantage of this choice is that it forces us to introduce the
notion of standard model, preparing a corollary of incompleteness: the
existence of non standard models. Another is that having a standard
model ${\cal M}$ permits, as we shall see, to prove easily the weak
representation theorem.

As for many theorems, we should not attempt to have the weakest
hypotheses in an elementary course, and the discussion of the most
general form of the theorem can be left for an advanced course.

\section{Which language ?}

The choice of the language is a delicate problem. To be able to
associate to each natural number $n$ a term $\underline{n}$, we need
to assume that the language in which the theory is expressed contains
enough function symbols to express the natural numbers. This forbids
to apply the incompleteness theorem to theories, such as set theory,
that has no function symbols.

An alternative is to use a theory where it is possible to construct
propositions with free variables characterizing the natural number
$0$, the successor relation, addition, multiplication, etc. such that
some propositions---essentially the axioms of Robinson's
arithmetic---are provable.  Then, for each natural number $p$, it is
possible to build a proposition $N_p$ characterizing the number $p$
and write $\fa x~(N_p[x] \Rightarrow A[x])$ or $\ex x~(N_p[x] \wedge
A[x])$ instead of $A[\underline{p}]$.

This is the choice made in \cite{Dowek}. In retrospect, it is a bit
heavy and it could have been left for an exercise.

\section{Numbering}

A first step in a proof of an incompleteness theorem is often the
introduction of the notion of numbering. For the students, this notion
is often both surprising and trivial. Surprising because theories,
algorithms, computers, etc. manipulate various datatypes without
numbering them. Trivial because, today, everyone knows
that texts, images, sounds, etc. are eventually coded as digits, hence
numbers.

Often, propositions, proofs, programs, Turing machines, etc.  are
numbered independently of each other. This leads to introduce many
definitions.  We defend that this notion of numbering should be made
general enough, so that the numbering of propositions, proofs,
etc. are just instances of this general definition.

On the other hand, this notion of numbering cannot be made too general
because composing a numbering with a non computable function yields
another function, that should not be considered as a numbering
\cite{BokerDershowitz}.

In \cite{Dowek}, we have proposed to restrict numberings to
articulated sets, using the notion, common in linguistics
\cite{Martinet}, of articulation. A set is said to be $0$-articulated
when it is finite. It is said to be $(n+1)$-articulated when it is a
set of finite trees labeled with elements of a $n$-articulated set.
So, the set containing the connectors, the quantifiers, the function
symbols, and the predicate symbols a of theory is $0$-articulated, the
set containing these symbols and the variables is 
$1$-articulated, the set of terms and propositions
is $2$-articulated, the set of sequents is $3$-articulated, the set of
proofs is $4$-articulated, etc. More generally, all the objects we
need to number are elements of an articulated set.

Using Cantor's bijection $;$ from ${\mathbb N}^2$ to ${\mathbb N}
\setminus \{0\}$
$$n;p = (n+p)(n+p+1)/2+ n + 1$$ we can number any tree of a
$(n+1)$-articulated set, $f(t_1, ..., t_p)$, whose root is labeled
with $f$ and whose immediate subtrees are $t_1$, ..., $t_p$ as
$$\quot{f(t_1, ..., t_p)} = \quot{f} ; (\quot{t_1} ; ( ...;
(\quot{t_p} ; 0)...))$$ where the first $\quot{.}$ is the numbering of
the $n$-articulated set the label $f$ belongs to and the others are
the function currently defined by induction. This way, the numbering
depends only on the numbering of the elements of the $0$-articulated
set we start with. And it is easy to prove that, as this set is
finite, changing this numbering does not change the set of computable
functions from an articulated set to another.

This is an example of an abstraction mechanism that permits to avoid
arbitrary definitions and tedious repetitions. Of course, as all
general notions, it must be illustrated with concrete examples, but
these concrete examples should not replace the general notion in the
definitions and the proofs.

Using this notion of articulation also forces us to define all the
objects we want to number---propositions, proofs, programs, etc.---as
trees. This means that the abstract syntax of these objects is
emphasized, and not their concrete syntax. We believe this is a good
thing and the notions of parenthese, unique reading, prefix, infix,
and postfix operator, precedence, etc., that are not specific to
logic languages should not be addressed in a logic course, but in a
language theory course.

Finally, when a proposition $A$ is numbered as $\quot{A}$, this number
$n$ is often used to build the term $S^n(0)$, written
$\underline{n}$. So, many expressions in the proofs are of the form
$\underline{\quot{A}}$.  Although the functions $\quot{.}$ and
$\underline{.}$ are sometimes used in isolation, introducing a
notation $\reflect{A}$ for $\underline{\quot{A}}$ clarifies the
proofs. Indeed, this composition of the functions $\quot{.}$ and
$\underline{.}$, mapping a proposition of arithmetic to terms of
arithmetic, is the homogeneous reflection notion.

\section{G\"odel's $\beta$ function and the definition of computable functions}

A second step in a proof of an incompleteness theorem is often the
association of a proposition to each computable function.

In fact, we do not associate one proposition to each computable
function, but to each construction of a computable function. For
instance, the binary null function and the composition of the unary
null function with a binary projection are extensionnally equal, but
different propositions are associated to these constructions.

We thus can introduce first a notion of construction of a computable
function, in such a way that $Z^2$ and $\circ^2_1(Z^1,\pi^2_1)$ are
different constructions of the same function. Such a construction is a
tree, labeled with symbols $Z^n$ for the zero function of arity $n$,
{\it Succ} for the successor function, $\mu^n$ for the minimization of
a function of arity $n+1$, etc. It may be called a ``program'', as it
is a syntactic object expressing a computable function.  It is in fact
the derivation tree, labeled with rule names \cite{Dowek} associated
to the inductive definition of computable functions, that is a proof
that the function is computable.

The set of computable functions is often defined as the smallest set
containing
\begin{itemize}
\item the projections, 
\item the null function,
\item and the successor function
\end{itemize}
and closed by 
\begin{itemize}
\item composition,
\item definitions by induction, 
\item and minimization,
\end{itemize}
leading to the language $\pi^n_i$, $Z^n$, {\it Succ}, $\circ^n_i$,
$\mbox{\it Rec}^n$, $\mu^n$.

Associating a proposition to each program is straightforward---$y = 0$
for the program $Z^n$, $y = S(x_1)$ for program $Succ$, etc.---except
for the definitions by induction, that require the use of G\"odel's
$\beta$ function and the Chinese remainder theorem.

This difficulty can be avoided, if we use an alternative definition of
the set of computable functions, adding three more functions:
addition, multiplication, and the characteristic function of the order
relation, that is the function $\chi_{\leq}$, such that
$\chi_{\leq}(n,p) = 1$ if $n \leq p$, and $\chi_{\leq}(n,p) = 0$
otherwise, and dropping definitions by induction.

The representation of programs is simplified.

\begin{itemize}
\item 
To $\pi_i^n$, we associate the proposition $y = x_i$.
\item 
To $Z^n$, we associate the proposition $y = 0$.
\item
To $\mbox{\it Succ}$, we associate the proposition $y = S(x_1)$.
\item
To $+$, we associate the proposition $y = x_1 + x_2$.
\item
To $\times$, we associate the proposition $y = x_1 \times x_2$.
\item
To $\chi_{\leq}$, we associate the proposition 
$$(x_1 \leq x_2 \wedge y = 1) \vee (x_2 < x_1 \wedge y = 0)$$ 
where the proposition $x \leq y$ abbreviates $\ex z~(z + x = y)$ and
$x < y$ abbreviates $S(x) \leq y$.

\item
To $\circ^n_m(h, g_1, \ldots, g_m)$, we associate the proposition 
$$\ex w_1 \ldots \ex w_m~(B_1[x_1, \ldots, x_n, w_1] \wedge \ldots
\wedge B_m[x_1, \ldots, x_n, w_m] \wedge C[w_1, \ldots, w_m, y])$$
where $B_1$, ..., $B_m$, and $C$, represent the programs $g_1$, ...,
$g_m$, and $h$.

\item
To $\mu^n(g)$, we associate the proposition 
$$\fa z~(z < y \Rightarrow \ex w~(B[x_1, \ldots, x_n, z, S(w)]))
\wedge B[x_1, \ldots, x_n, y, 0]$$ where $B$ represents the program
$g$.
\end{itemize}
Fundamentally, what makes this notion of representation easy is the
similarity between the symbols allowing to construct terms in
arithmetic: variables, $0$, $S$, $+$, and $\times$, and five of the
eight clauses defining the set of computable functions: projections,
the null functions, the successor function, addition, and
multiplication.

But, of course, we need to prove the equivalence of these two
definitions of the set of computable functions, and this requires the
use of G\"odel's $\beta$ function and the Chinese remainder theorem.
Proving this equivalence \cite{Dowek}\footnote{The original proof
  contained a few gaps, a corrected proof---in French---is available
  online\\ {\tt
    http://www.lsv.fr/$\tilde{~}$dowek/Books/Lc/prop317.pdf}} is even
slightly more difficult than directly using G\"odel's $\beta$ function
and the Chinese remainder theorem to represent programs as
propositions, but it makes the proof more modular.  When proving this
equivalence, only functions mapping natural numbers to natural numbers
are used and the notion of proposition is not mentioned. This
equivalence can also be proved long before the incompleteness theorem
is discussed.  It can be motivated by other goals than the
incompleteness theorems: for instance it simplifies the proof of other
theorems such as the representation theorem of computable functions as
rewrite systems, terms of the $\lambda$-calculus \cite{Krivine},
Turing machines, etc.  When computability is taught in a different
course than logic, this equivalence should, of course, be taught in
the computability course, and not in the logic course.

Computable functions could even be defined with these eight clauses,
and G\"odel's $\beta$ function and the Chinese remainder theorem would
then be used only to prove that the set of computable functions is
closed by definitions by induction.

\section{The form of the representation theorem}

\subsection{The weak representation theorem}
\label{weakrep}

Let ${\cal T}$ be a theory that has a decidable set of axioms, is an
extension of Robinson's arithmetic, and has a standard model ${\cal
  M}$.

We want to prove that if a proposition $A$ represents a program $f$,
then the proposition $A[\underline{p_1}, ..., \underline{p_n},
  \underline{q}]$ expresses that the program $f$ terminates at $p_1$,
..., $p_n$ and returns $q$, that is that this proposition is provable
in ${\cal T}$ if and only if $f(p_1, ..., p_n) = q$.

A simple induction on the structure of the program permits to prove
that if $f(p_1, ..., p_n) = q$, then the proposition
$A[\underline{p_1}, ..., \underline{p_n},\underline{q}]$ is provable
in ${\cal T}$.  The completeness theorem shows that if this
proposition is provable in ${\cal T}$, it is valid in ${\cal M}$. And
a simple induction on the structure of the program $f$ shows that, if
this proposition is valid in ${\cal M}$, then $q = f(p_1, ..., p_n)$.

\subsection{The undecidability of provability}

In the same way, we can prove that the proposition $\ex
y~A[\underline{p_1}, ..., \underline{p_n},y]$ expresses that the
program $f$ terminates at $p_1$, ..., $p_n$, that is that this
proposition is provable in ${\cal T}$ if and only if $f$ terminates at
$p_1, ..., p_n$.

The computable function $F$ mapping $f$ and $p_1, ..., p_n$ to the
proposition $\ex y~A[\underline{p_1}, ..., \underline{p_n},y]$ thus
reduces the halting problem to provability in ${\cal T}$.

If $G$ is the function mapping a proposition $A$ to $1$ if it is
provable and to $0$ otherwise, then the function $G \circ F$ maps $f$
and $p_1, ..., p_n$ to $1$ is $f$ terminates at $p_1, ..., p_n$, and
to $0$ otherwise.  Using the contrapositive of the closure of
computable functions by composition, as $G \circ F$ is not computable
and $F$ is, $G$ is not. Thus, the undecidability of provability in
${\cal T}$ is a mere consequence of the undecidability of the halting
problem and of this representation theorem.

\subsection{A stronger representation theorem}

In some proofs of the incompleteness theorems, we need a stronger
theorem expressing that if $f$ terminates at $p_1, ..., p_n$ then the
proposition
$$\fa y~(A[\underline{p_1}, ..., \underline{p_n},y] \Leftrightarrow y
= \underline{f(p_1,...,p_n)})$$ 
is provable. For this theorem, we do not need the theory ${\cal T}$ to
have a standard model, or even to be consistent.  All we need it that
it has a decidable set of axioms an it is an extension of Robinson's
arithmetic.

Note that this representation theorem is both stronger and weaker than
that of Section \ref{weakrep}. It is stronger because the equivalence
with $f$ is internalized, it is expressed by an equivalence in the
language, that holds for all $y$. But, it is weaker because it says
nothing when $f$ does not terminate at $p_1$, ..., $p_n$, while the
theorem of Section \ref{weakrep} shows that $A[\underline{p_1}, ...,
  \underline{p_n},\underline{q}]$ is not provable in this case.

But, assuming that the theory ${\cal T}$ is moreover consistent, when
the function terminates at $p_1, ..., p_n$, the weak representation
theorem is a consequence of the strong one as
$$A[\underline{p_1}, ..., \underline{p_n},\underline{q}]$$
is provable if and only if $\underline{q} = \underline{f(p_1, ..., p_n)}$ 
is, that is if and only if $q = f(p_1, ..., p_n)$.

The statement of this strong representation theorem is slightly less
natural than that of the weak one, as it introduces a asymmetry
between the arguments and the value of the program. The arguments are
closed terms $\underline{p_1}$, ..., $\underline{p_n}$ quantified
outside the language, while the value is a variable $y$ quantified in
the language itself.

As a consequence, if the proof of the strong representation is direct
for the seven of the eight cases, it is slightly less direct for
minimization, $\mu^n(g)$, where an argument of the function $g$
becomes the value of the function $\mu^n(g)$.  We need to prove
$$\fa y~((\fa z~(z < y \Rightarrow \ex w~B[\underline{p_1}, \ldots,
  \underline{p_{n}}, z, S(w)]) \wedge B[\underline{p_1}, \ldots,
  \underline{p_n}, y, 0]) \Leftrightarrow y = \underline{r})$$
assuming $\mu^n(g)$ terminates and takes the value $r$ at $p_1$, ...,
$p_n$.  Proving
$$\fa y~(y = \underline{r} \Rightarrow (\fa z~(z < y \Rightarrow \ex
w~B[\underline{p_1}, \ldots, \underline{p_{n}}, z, S(w)]) \wedge
B[\underline{p_1}, \ldots, \underline{p_n}, y, 0]))$$ 
that is equivalent to
$$\fa z~(z < \underline{r} \Rightarrow \ex w~B[\underline{p_1},
\ldots, \underline{p_{n}}, z, S(w)]) \wedge B[\underline{p_1},
\ldots, \underline{p_n}, \underline{r}, 0]$$ 
is easy, as the bounded quantification can be reduced to a finite
conjunction. But to prove the converse
$$\fa y~( (\fa z~(z < y \Rightarrow \ex w~B[\underline{p_1}, \ldots,
  \underline{p_{n}}, z, S(w)]) \wedge B[\underline{p_1}, \ldots,
  \underline{p_n}, y, 0]) \Rightarrow y = \underline{r})$$ 
we need to use the fact that from the hypothesis
$\fa z~(z < y \Rightarrow \ex
w~B[\underline{p_1}, \ldots, \underline{p_{n}}, z, S(w)])$ we can
deduce 
$\underline{r} < y \Rightarrow \ex w~B[\underline{p_1}, \ldots, 
\underline{p_{n}}, \underline{r}, S(w)]$ 
thus 
$\underline{r} < y \Rightarrow \ex w~0 = S(w)$ 
and 
$y \leq \underline{r}$, to show that it is sufficient to prove the proposition
$$\fa y~(y \leq \underline{r} \Rightarrow (B[\underline{p_1}, \ldots, 
  \underline{p_n}, y, 0] \Rightarrow y = \underline{r}))$$
that reduces to a finite conjunction.

As we shall see, the weak representation theorem is enough for several
proofs of the first incompleteness theorem and in an elementary
course, we can restrict to this theorem, while the strong theorem is
needed for an advanced course.

\section{The various proofs of the first incompleteness theorem}

There are two families of proofs of the first incompleteness theorem. 

\subsection{The computer scientist's proofs}

In the first family, the incompleteness of the theory is seen as a
consequence of the undecidability of provability in this theory.  Let
${\cal T}$ be a theory that has a decidable set of axioms, is an
extension of Robinson's arithmetic, and has a standard model ${\cal
M}$.  As we have seen, provability in ${\cal T}$ is undecidable.

If this theory were complete, then the computable function mapping
$\quot{A}$ to the least $x$ such that $\mbox{\it proof}(x,\quot{A}) =
1$ or $\mbox{\it proof}(x,\quot{\neg A }) = 1$, where the computable
function {\it proof} maps $n$ and $p$ to $1$ if $n = \quot{\pi}$, $p =
\quot{A}$ and $\pi$ is a proof of $A$, and to $0$ otherwise, would be
total and would return a proof of $A$ if and only if $A$ is
provable. Thus it would permit to build an algorithm deciding
provability in ${\cal T}$.

This theorem can even be made more abstract as the fact that a
semi-decidable set whose complement is also semi-decidable is
decidable.

This proof uses many notions of theoretical computer science: the notion of
computable function, the notion of reduction, and the notion of proof
search: the function mapping $\quot{A}$ to the least $x$ such that
$\mbox{\it proof}(x,\quot{A}) = 1$ is a generate-and-test
proof search algorithm and that mapping $\quot{A}$ to the least $x$
such that $\mbox{\it proof}(x,\quot{ A}) = 1$ or $\mbox{\it
  proof}(x,\quot{\neg A}) = 1$ is a similar algorithm searching
simultaneously for a proof of $A$ and $\neg A$. So, this proof may be
called the computer scientist's proof.

Considering a theory that has a decidable set of axioms, is an
extension of Robinson's arithmetic, or Peano arithmetic, and has a
standard model, proving the weak representation theorem for this
theory, deducing the undecidability of provability in this theory, and
then its incompleteness as a corollary is probably a sufficient goal
for an elementary course, giving a first exposition to incompleteness.

\subsection{The quick mathematician's proof}

The second family of proofs effectively constructs a proposition $G$
such that neither $G$ nor $\neg G$ are provable.

It gives a less central {\it r\^ole} to the notion of computable
function. Often the representation of a few functions such as the
function {\it proof}, the substitution function mapping $n$ and $p$ to
$m$ if $n = \quot{A}$ and $m = \quot{(\underline{p}/x)A}$, and the
negation function mapping $n$ to $m$ if $n = \quot{A}$ and $m =
\quot{\neg A}$ are sufficient.

A simple formulation of this proof, which can be given as an exercise,
even in an elementary course, is to consider only one function, that
is a mixture of the {\it proof} function and of the substitution
function, mapping $n$, $p$, and $q$ to $1$ if $n = \quot{\pi}$, $p =
\quot{A}$ and $\pi$ is a proof of $(\underline{q}/x)A$.

Let ${\cal T}$ be a theory that has a decidable set of axioms, is an
extension of Robinson's arithmetic, and has a standard model ${\cal
  M}$.

Calling $F$ the proposition representing the function $f$, $T$ the
diagonal proposition \linebreak $\fa x~\neg F[x,w,w,\underline{1}]$
and $G$ the proposition $T[\reflect{T}]$, it is easy to prove that
neither $G$ nor $\neg G$ is provable, in a theory that has a decidable
set of axioms, is an extension of Robinson's arithmetic, and has a
standard model.

If $G$ is provable, then $\fa x~\neg
F[x,\reflect{T},\reflect{T},\underline{1}]$ is provable, thus for all
$n$, $f(n,\quot{T},\quot{T}) = 0$, thus $T[\reflect{T}]$, that is $G$,
is not provable. A contradiction.  If $\neg G$ is provable, then there
exists a natural number $n$ such that $f(n,\quot{T},\quot{T}) = 1$
thus $G$ is provable, thus $\bot$ is provable. A contradiction.

\subsection{The thorough mathematician's proof}

This quick proof is still a bit mysterious because the {\it proof}
function and the substitution function are mixed in this function $f$
and the meaning of the proposition $G$---its relation to the liar's
paradox, to the diagonal argument, and to self reference---are not
explicit.

A more thorough proof, which also prepares the proof of the second
incompleteness theorem better, is to use the {\it proof} function and
the substitution function $s$ to decompose this function $f$.

Then, we can introduce a provability proposition 
$$\mbox{\it Bew} = \ex z~\mbox{\it Proof}[z,x,\underline{1}]$$ and
introduce the notation, inspired by provability logic, $\Box A$ for
$\mbox{\it Bew}[\reflect{A}]$.

It is easy to prove the necessitation lemma: if $A$ is provable, then
$\Box A$ is provable.  Indeed, if $A$ has a proof $\pi$, then, by the
representation theorem, $\mbox{\it
  Proof}[\reflect{\pi},\reflect{A},\underline{1}]$. Thus, $\ex
z~\mbox{\it Proof}[z,\reflect{A},\underline{1}]$, that is $\mbox{\it
  Bew}[\reflect{A}]$, that is $\Box A$, is provable.

The converse of this lemma is false in general: a theory can prove
$\Box \bot$ without proving $\bot$, but this converse holds if the theory is
$\omega$-consistent. Indeed from $\Box A$, that is $\ex z~\mbox{\it
  Proof}[z,\reflect{A},\underline{1}]$ we can deduce that there exists
an $n$ such that $\neg \mbox{\it
  Proof}[\underline{n},\reflect{A},\underline{1}]$ is not provable,
thus using the strong representation theorem $\neg \underline{1} =
\mbox{\it proof}(n,\quot{A})$ is not provable, thus $\mbox{\it
  proof}(n,\quot{A}) \neq 0$, $\mbox{\it proof}(n,\quot{A}) = 1$ and
$A$ is provable.

Let ${\cal T}$ be a theory that has a decidable set of axioms, is an
extension of Robinson's arithmetic, and is $\omega$-consistent.

Using the strong representation theorem, we can now build the liar's 
proposition $G$ such that $G \Leftrightarrow \neg \Box G$ is provable. 
To do so, we introduce a proposition $D$ expressing that 
$A[\underline{p}]$ is not provable, 
$$\ex z~(\neg \mbox{\it Bew}[z] \wedge S[x_1, x_2,z])$$
a proposition $E$ expressing
that $A[\reflect{A}]$ is not provable 
$$D[w,w]$$
and $G$ 
$$E[\reflect{E}]$$
that is 
$$\ex z~(\neg \mbox{\it Bew}[z] \wedge S[\reflect{E},\reflect{E},z])$$
then $G$ is provably equivalent to 
$$\ex z~(\neg \mbox{\it Bew}[z] \wedge z = \reflect{E [\reflect{E}]})$$
that is 
$$\neg \mbox{\it Bew}[\reflect{G}]$$

So, the proposition $G$ is a fixed point of the proposition
$\neg \mbox{\it Bew}[x]$. 
This theorem can be generalized \cite{Miquel}. 
We consider 
any proposition $C$ containing a variable $x$ and we 
prove that there exists 
a proposition $G$ such that $G \Leftrightarrow C[\reflect{G}]$
is provable. Instead of defining $D$ as 
$\ex z~(\neg \mbox{\it Bew}[z] \wedge S[x_1,x_2,z])$
we define it as 
$$\ex z~(C[z] \wedge S[x_1,x_2,z])$$
$E$ as 
$$D[w,w]$$
and $G$ as 
$$E[\reflect{E}]$$
that is 
$$\ex z~(C[z] \wedge S[\reflect{E},\reflect{E},z])$$
The proposition $G$ is provably equivalent to 
$$\ex z~(C[z] \wedge z = \reflect{E [\reflect{E}]})$$
hence to 
$$C[\reflect{G}]$$

Taking the proposition $C = \neg \mbox{\it Bew}[x]$, 
we get the liar's proposition above.

Taking the proposition $C = \mbox{\it Bew}[x]$, we get Henkin's
truth-teller's proposition $H$ such that $H \Leftrightarrow \Box H$ is
provable.

Taking the proposition $C = \mbox{\it Bew}[x] \Rightarrow P$, we
get L\"ob's proposition such that $L \Leftrightarrow (\Box L
\Rightarrow P)$ is provable.

Taking the proposition
$$C = \fa y~(\mbox{\it Proof}[y,x] \Rightarrow \ex z~(z \leq y
\wedge (\ex w~(Neg[x, w] \wedge \mbox{\it Proof}[z,w]))))$$ 
we get Rosser's proposition such that
$$R \Leftrightarrow \fa y~(\mbox{\it Proof}[y,\reflect{R}] \Rightarrow
\ex z~(z \leq y \wedge (\ex w~(Neg[\reflect{R},w] \wedge \mbox{\it
  Proof}[z,w]))))$$ is provable,
that is 
$$R \Leftrightarrow \fa y~(\mbox{\it Proof}[y,\reflect{R}] \Rightarrow
\ex z~(z \leq y \wedge (\mbox{\it Proof}[z,\reflect{\neg R}])))$$ 
is provable.

This fixed point theorem is another example where
abstraction---considering an arbitrary proposition---simplifies proofs,
avoiding redundancy.

The proof that neither $G$ nor $\neg G$ are provable in ${\cal T}$ is
then quite direct using $G \Leftrightarrow \neg \Box G$, the
necessitation and its converse: if $G$ is provable then $\Box G$
(necessitation) and $\neg \Box G$ (equivalence) also. Thus $\bot$ is
provable. A contradiction.  If $\neg G$ is provable, then $\Box G$
also (equivalence) hence $G$ also (converse of necessitation). Thus
$\bot$ is provable. A contradiction.

As this proof uses the converse of necessitation, it requires the
theory ${\cal T}$ to be $\omega$-consistent.

\section{Consistency and $\omega$-consistency}

The three proofs discussed above use an hypothesis stronger than
consistency: $\omega$-consistency or the existence of a standard
model.  The G\"odel-Rosser theorem permits to weaken this hypothesis
to consistency. Thus, we consider a theory ${\cal T}$ that has a
decidable set of axioms, is an extension of Robinson's arithmetic, and
is consistent.

This proof uses ``the little converse of necessitation'': if the proposition 
$$\ex x~(x \leq \underline{n} \wedge \mbox{\it Proof}[x,\reflect{A}])$$
is provable then so is $A$. 
Indeed from 
$$\ex x~(x \leq \underline{n} \wedge \mbox{\it Proof}[x,\reflect{A}])$$
we can deduce 
$$\mbox{\it Proof}[\underline{0}, \reflect{A}] \vee \mbox{\it
  Proof}[\underline{1}, \reflect{A}] \vee ...  \vee \mbox{\it
  Proof}[\underline{n}, \reflect{A}]$$ and if $A$ were not be
provable, then $\neg \mbox{\it Proof}[\underline{0}, \reflect{A}]$,
$\neg \mbox{\it Proof}[\underline{1}, \reflect{A}]$, ..., $\neg
\mbox{\it Proof}[\underline{n}, \reflect{A}]$ would be provable. Then
$\bot$ would be provable and $A$ also.  A contradiction.

Then, the proof that neither $R$ not $\neg R$ is provable is quite direct.
If $R$ has a proof $\pi$ then the proposition
$$\mbox{\it Proof}[\reflect{\pi},\reflect{R}]$$
is provable and the proposition 
$$\fa y~(\mbox{\it Proof}[y, \reflect{R}] \Rightarrow
\ex z~(z \leq y \wedge \mbox{\it Proof}[z, \reflect{\neg R}]))$$
also. Thus, the proposition 
$$\ex z~(z \leq \reflect{\pi} \wedge \mbox{\it Proof}[z, \reflect{\neg R}])$$
is provable and, by the little converse of necessitation, $\neg R$ also. 
Thus $\bot$ is provable, contradicting the consistency of the theory. 

If $\neg R$ has a proof $\pi$, then the proposition 
$$\mbox{\it Proof}[\reflect{\pi}, \reflect{\neg R}]$$
is provable and the proposition 
$$\ex y~(\mbox{\it Proof}[y, \reflect{R}] \wedge
\fa z~(\neg z \leq y \vee \neg \mbox{\it Proof}[z, \reflect{\neg R}]))$$
also. So, the proposition 
$$\ex y~(\mbox{\it Proof}[y,\reflect{R}] \wedge
(\neg \reflect{\pi} \leq y \vee \neg \mbox{\it Proof}
[\reflect{\pi},\reflect{\neg R}]))$$
is provable. Hence the proposition 
$$\ex y~(y < \reflect{\pi} \wedge \mbox{\it Proof}[y, \reflect{R}])$$
is provable and by the little converse of necessitation $R$ also.
Thus $\bot$ is provable, contradicting the consistency of the theory.

Introducing this little converse of necessitation permits to 
make this proof more modular. 

Yet, because, not only the notion of provability, but also the notion
of proof is used in this proof, it cannot be formulated abstractly
using the notation $\Box A$.

\section{The second incompleteness theorem}

The second incompleteness theorem shows that a consistent extension of
Peano arithmetic does not prove its own consistency.

\subsection{The Hilbert-Bernays lemmas} 

Besides the necessitation lemma: if $A$ is provable, so is $\Box A$,
the proof of the second incompleteness theorem requires two more
lemmas: the internalization of {\it modus ponens}
$$\Box (A \Rightarrow B) \Rightarrow \Box A \Rightarrow \Box B$$
is provable and the internalization of necessitation 
$$\Box A \Rightarrow \Box \Box A$$
is provable.

The proofs of these two lemmas use induction and cannot be proved in
Robinson's arithmetic.  So, we consider a theory ${\cal T}$ that has a
decidable set of axioms, is an extension of Peano arithmetic, and is
consistent.

Many books just say that these two lemmas have long and tedious
proofs.  Some insist on the fact the informal statements
\begin{center}
If $A \Rightarrow B$ and $A$ are provable, then $B$ is provable.
\end{center}
and 
\begin{center}
If $A$ is provable, then $\Box A$ is provable.
\end{center}
can be expressed in arithmetic as the two propositions above and that
their proofs can be formalized in arithmetic---formalizing these
proofs being long and tedious.

As remarked by Miquel \cite{Miquel}, these proofs depend on the choice
of the proposition {\it Proof}, that is of the proof-checking program
expressing the function {\it proof}. The program we use transforms
trees whose leaves are labeled with proofs or with the symbols $0$ and
$1$ and whose internal nodes are all labeled with a symbol {\it and}.

In such a tree we consider two kinds of reducible expressions.
\begin{itemize}
\item An internal node labeled with {\it and} whose both children are
  labeled with $0$ or $1$ is a reducible expression. It reduces
  to $1$,
  when both children are labeled with $1$, and to $0$ otherwise.
\item A leaf labeled with a proof $\pi$ is a reducible expression.  
If this proof has
  a root labeled with the proposition $A$ and immediate subproofs
  $\pi_1$, ..., $\pi_n$, we write it $A(\pi_1,...,\pi_n)$.  Let $B_1$,
  ..., $B_n$ be the propositions labelling the roots of $\pi_1$, ...,
  $\pi_n$.  If there is a deduction rule allowing to deduce $A$ from
  $B_1$, ..., $B_n$, then it reduces to the tree
  $(...(\pi_1~\mbox{\it and}~\pi_2) ...~\mbox{\it and}~\pi_n)$ and
  to $1$ if $n = 0$. Otherwise, it reduces to $0$.
\end{itemize}

We first define a program {\it step} that reduces the leftmost reducible 
expression
in a tree. We then define a program {\it check} iterating the program
{\it step} until obtaining an irreducible tree, containing just one
node, labeled with $1$ or with $0$, which is the result of the
algorithm.

The program {\it proof} applied to $\pi$ and $A$ applies the program
{\it check} to the tree containing just one node labeled with
$\pi$. If the result is $1$ and the root of $\pi$ is $A$, it returns
$1$, otherwise it returns $0$.

If it easy to prove, by induction on $n$ that if $\mbox{\it
  step}^n(\quot{\pi_1}) = \quot{\pi'_1}$ then $\mbox{\it
  step}^n(\quot{\pi_1~\mbox{\it and}~\pi_2}) = \quot{\pi'_1~\mbox{\it
    and}~\pi_2}$ and that if $\mbox{\it step}^n(\quot{\pi_2}) =
\quot{\pi'_2}$ then $\mbox{\it step}^n(\quot{1~\mbox{\it and}~\pi_2} =
\quot{1~\mbox{\it and}~\pi'_2}$.

Thus, if $\mbox{\it proof}(\quot{\pi_1},\quot{A \Rightarrow B}) = 1$
and $\mbox{\it proof}(\quot{\pi_2},\quot{A}) = 1$ then there exists
natural numbers $n$ and $p$ such that $\mbox{\it step}^n(\quot{\pi_1})
= \quot{1}$ and $\mbox{\it step}^p(\quot{\pi_2}) = \quot{1}$.  Thus,
$\mbox{\it step}^{n+p}(\quot{\pi_1~\mbox{\it and}~\pi_2}) =
\quot{1~\mbox{and}~1}$ and $\mbox{\it
  step}^{n+p+1}(\quot{\pi_1~\mbox{\it and}~\pi_2}) = \quot{1}$.  As
$\mbox{\it step}(\quot{B(\pi_1,\pi_2)} = \quot{\pi_1~\mbox{\it
    and}~\pi_2}$, \linebreak $\mbox{\it
  step}^{n+p+2}(\quot{B(\pi_1,\pi_2)}) = \quot{1}$.  Thus
$proof(\quot{B(\pi_1,\pi_2)},\quot{B}) = 1$.

This proof, which uses induction only, can be expressed in Peano
arithmetic and it is a proof of
$$\fa x_1 \fa x_2 \fa y~(\mbox{\it Proof}[x_1,\reflect{A \Rightarrow
B}, \underline{1}] \Rightarrow \mbox{\it Proof}[x_2,\reflect
A,\underline{1}] \Rightarrow M[\reflect{B},x_1,x_2,y] \Rightarrow
\mbox{\it Proof}[y,\reflect{B},\underline{1}])$$ where $M$ is the
propositions representing the function mapping $\quot{B}$,
$\quot{\pi_1}$ and $\quot{\pi_2}$ to $\quot{B(\pi_1,\pi_2)}$.

From this proposition and the totality of the function mapping
$\quot{B}$, $\quot{\pi_1}$ and $\quot{\pi_2}$ to
$\quot{B(\pi_1,\pi_2)}$, we get a proof of
$$(\ex x_1 \mbox{\it Proof}[x_1,\reflect{A \Rightarrow B},\underline{1}])
\Rightarrow (\ex x_2 \mbox{\it Proof}[x_2,\reflect A,\underline{1}])
\Rightarrow  \ex y~\mbox{\it Proof}[y,\reflect{B},\underline{1}]$$
that is of
$$\Box(A \Rightarrow B) \Rightarrow \Box A \Rightarrow  \Box B$$
which is the internalization of the {\it modus ponens}. 

A similar argument can be given for all the other deduction rules.
This permits to give a new proof of the necessitation lemma: if $A$
has a proof, so does $\Box A$, building the proof of $\Box A$, step by
step, by induction on the proof of $A$.

This proof also can be formalized in arithmetic as a proof of
$$\Box A \Rightarrow \Box \Box A$$

\subsection{The second incompleteness theorem}

The proof of the second incompleteness theorem: the theory ${\cal T}$
does not prove its own consistency $\neg \Box \bot$, is then quite
direct.

Let $G$ a the proposition such that $G \Leftrightarrow \neg \Box G$ is
provable.

The proposition $G \Rightarrow \Box G \Rightarrow \bot$ is provable,
so, using the necessitation lemma and the internalization of the {\it
  modus ponens}, the propositions $\Box(G \Rightarrow \Box G
\Rightarrow \bot)$ and $\Box G \Rightarrow \Box \Box G \Rightarrow
\Box \bot$ also. Thus, using the internalization of necessitation, the
proposition $\Box G \Rightarrow \Box \bot$ also.

It can be noticed that this proposition is the internalization of the
first half of the first incompleteness theorem: if $G$ is provable,
then $\bot$ also.

Now, if we assume that the proposition $\neg \Box \bot$ is provable,
we can deduce that the proposition $\neg \Box G$ is also
provable. Thus, the proposition $G$ also, and by necessitation, $\Box
G$ also, thus $\bot$ also, contradicting the consistency of the
theory.

\subsection{Provability logic} 

The proof of the first incompletenes theorem uses the symbol 
$\Box$ but only the necessitation lemma. The proof of the second 
theorem, in contrast, uses more modal logic: the lemmas $K$ and $4$
$$\Box (A \Rightarrow B) \Rightarrow \Box A \Rightarrow \Box B$$
$$\Box A \Rightarrow \Box \Box A$$ 

It is possible to mention modal logic here, and in particular
provability logic, but the notation of provability logic and its
modularity---proving first the Hilbert-Bernays lemmas, and then using
them in the proof of the second incompleteness theorem---can be used
without defining provability logic {\it per se}.

\subsection{L\"ob's theorem} 

Using the notations of provability logic and a general fixed point
theorem, L\"ob's theorem is a straightforward extension of the second
incompleteness theorem.  L\"ob's theorem is: for any proposition $P$,
if $\Box P \Rightarrow P$ is provable then $P$ is.

Instead of the proposition $G$, we use the proposition $L$ such that
$L \Leftrightarrow (\Box L \Rightarrow P)$ is provable.

The proposition $L \Rightarrow \Box L \Rightarrow P$ is provable, so,
using the necessitation lemma and the internalization of the {\it
  modus ponens}, the propositions $\Box(L \Rightarrow \Box L
\Rightarrow P)$ and $\Box L \Rightarrow \Box \Box L \Rightarrow \Box
P$ also. Using the internalization of necessitation, the proposition
$\Box L \Rightarrow \Box P$ also.

Now, if we assume that the proposition $\Box P \Rightarrow P$ is
provable, we can deduce that the proposition $\Box L \Rightarrow P$ is
provable. Thus, the proposition $L$ also, and by necessitation $\Box
L$ also, thus $P$ is provable.

The second incompleteness theorem is a corollary of this theorem
taking $P = \bot$, so, this theorem can also be proved before the
second incompleteness theorem.

Finally, taking $P = H$, Henkin's truth-teller's proposition, such
that $H \Leftrightarrow \Box H$ is provable, another corollary is that
this proposition $H$ is provable. Indeed, as $\Box H \Rightarrow H$ is
provable, so is $H$.

\section{History and philosophy}

As the subject itself, the history of the incompleteness theorems is
very rich. One thing the students can learn is that the notion of
computable function, and the undecidability of provability, and hence
the computer scientist's proof came in the work of Church and Turing
(1936) after the incompleness theorems (1931). So, historically, the
first proof is the mathematician's proof. Another thing the students
can learn is the problems the second incompleteness theorem and the
undecidability solved: Hilbert's second problem and Hilbert's {\it
  Entscheidungsproblem} respectively.

On the more philosophical side, many commentators see in the
incompleteness a proof that a theory cannot speak about itself. On the
contrary, the second incompleteness theorem exists because the
consistency of a theory can be formulated in the theory itself. The
incompleteness is also often presented as a disaster. If it shows some
limits to the deductive method, it is not the end of it.

\section{Conclusion}

As suggested in the Introduction, the incompleteness theorems are rich
subject that cannot be taught in one course. We have tried to separate
what can be taught in an elementary course: the weak representation of
computable functions, the undecidability of provability, the first
incompleteness theorem, under reasonable hypotheses, from what can be
kept for an advanced course: the strong representation of computable
functions, the thorough proof of the first incompleteness theorem,
Hilbert-Bernays lemmas, the second incompleteness theorem, L\"ob's
theorem, and the minimal hypotheses to be used in these theorems.

We have also defended that the notion of computable function, the
equivalence of its two definitions, using G\"odel's $\beta$ function
and the Chinese remainder theorem, and concrete syntax should be
taught independently and before the undecidability and incompleteness
of arithmetic.

Finally, we have defended an abstract and modular approach to these
proofs, using abstract syntax, articulation, universal numbering,
provability logic notations, the converse of necessitation, the little
converse of necessitation, a general fixed point lemma, and
Hilbert-Bernays lemmas. This work needs to be continued. In particular
the G\"odel-Rosser theorem and the proofs of the Hilbert-Bernays
lemmas are not yet abstract enough.

Abstract definitions and general lemmas should of course be
illustrated with concrete examples, but these concrete examples should
not replace them.

\nocite{CoriLascar}
\nocite{Andrews}

\bibliographystyle{plain}
\bibliography{teachinggodel}
\end{document}